\documentclass{amsproc}
\usepackage{graphicx,pstricks,amssymb}
\usepackage{mathrsfs}
\newtheorem{thm}{Theorem}[section]
\newtheorem{cor}[thm]{Corollary}
\newtheorem{lem}[thm]{Lemma}
\newtheorem{prop}[thm]{Proposition}
\theoremstyle{definition}

\theoremstyle{remark}

\numberwithin{equation}{section}
\newcommand{\commentout}[1]{}

\newcommand{\Tr}{\mathrm{Tr}}

\begin{document}

\title{Differential equations for deformed Laguerre polynomials}
\author{Peter J. Forrester and Christopher M. Ormerod}

\begin{abstract}
The distribution function for the first eigenvalue spacing in the Laguerre unitary ensemble of finite 
size may be expressed in terms of a solution of the fifth Painlev\'e transcendent. The generating 
function of a certain discontinuous linear statistic of the Laguerre unitary ensemble can similarly
be expressed in terms of  
a solution of the fifth Painlev\'e equation. The methodology used to derive 
these results rely on two theories regarding differential equations for orthogonal polynomial
systems, one involving isomonodromic deformations and the other ladder operators.
We compare the two theories by showing how either can be used to obtain a characterization
of a more general Laguerre unitary ensemble average in terms of the Hamiltonian system for
Painlev\'e V.
\end{abstract}

\maketitle

\section{Introduction}
\subsection{Objective}
Spacing probabilities and moments of characteristic polynomials for random matrix ensembles
with unitary symmetry are intimately related to semi-classical orthogonal polynomials. By developing the
theory of particular semi-classical polynomials, it has been possible to characterize these
random matrix quantities in terms of both discrete and continuous Painlev\'e equations
\cite{PS90a,Hi96a,IW01,BR00,BR01b,BD00,PS03,FW03a,FW06,Forrester,CF06,BasorChen, TW94c, AV95}. Two methods of developing
the theory for this purpose have emerged. One has been to use a formulation in terms of
Lax pairs for isomondromic deformations of linear differential equations \cite{Magnus1992}. The other has proceeded
via a theory of ladder operators for orthogonal polynomial systems 
\cite{IsmailChenLadder,Chen:discontinuous}. It is the
purpose of this paper to compare these two methods as they apply to the particular discontinuous
semi-classical weight
\begin{equation}\label{eq1:weight}
w(x) = (1 - \zeta \theta(x - t) ) (x - t)^\alpha x^\mu e^{-x},
\end{equation}
with support on $\mathbb R^+$, $\zeta < 1$, where $\theta(y)$ denotes the Heaviside function
$\theta(y) = 1$ for $y > 0$, $\theta(y) = 0$ otherwise.

\subsection{A random matrix context}
The weight \eqref{eq1:weight} is relevant to the Laguerre unitary ensemble LUE${}_N^\alpha$ specified by the
eigenvalue probability density function
\begin{equation}\label{eq1:product}
\frac{1}{C}\prod_{l=1}^N x_l^\alpha e^{- x_l} \prod_{1 \le j < k \le N} (x_k - x_j)^2, \, x_l > 0. 
\end{equation}
For $\alpha$ a non-negative integer, this is realized by so called Wishart matrices
$X^\dagger X$ where $X$ is an $n \times N$ $(n \ge N, \, \alpha = n - N)$ complex
matrix of independent standard complex Gaussian matrices (see e.g.~\cite{Fo02}). To see how
\eqref{eq1:product} relates to \eqref{eq1:weight}, consider the random matrix average
\begin{equation}\label{bc.3}
\Big \langle \prod_{l=1}^N (1 - \zeta \theta(x_l - t) )
|x_l - t|^\mu x_l^\alpha e^{-x_l} \Big \rangle_{{\rm LUE}_N^\alpha}.
\end{equation}
This average is proportional to the multiple integral
\begin{equation}\label{eq1:Delta}
\Delta_N = \frac{1}{N!} \int_0^\infty dx_1 \cdots \int_0^\infty dx_N \,
\prod_{l=1}^N w(x_l) \prod_{1 \le j < k \le N} (x_k - x_j)^2,
\end{equation}
where $w(x)$ is given by \eqref{eq1:weight}. Introducing the moments
\begin{equation}\label{eq0:moment}
\mu_n := \int_0^\infty w(x) x^n \, dx
\end{equation}
it is an easy result that
\begin{equation}\label{eq0:deltadetdef}
\Delta_N = \det [ \mu_{j+k-2} ]_{j,k=1,\dots,N}.
\end{equation}
Moreover, the $\{\Delta_n\}$ can be calculated through a recurrence linking the
orthogonal polynomials $\{p_n(x)\}$ associated with $w(x)$.

Introducing the bilinear form
\begin{equation}\label{eq0:bilinear}
\left< f,g \right> = \int_I f(x)g(x)w(x)dx
\end{equation}
the orthogonal polynomials are specified by the requirements that $p_n$ be a
polynomial of degree $n$ and
\begin{equation}\label{eq0:orthonormality}
\left< p_i,p_j \right> = \delta_{i,j}.
\end{equation}
By way of application of the Gram-Schmidt process, the weight function and its associated support completely specify the
coefficients of each polynomial. Furthermore the orthonormality of the polynomials implies
the three term recurrence relation \cite{Szego}
\begin{equation}\label{eq0:3termrecurrence}
a_{n+1}p_{n+1} = (x-b_n)p_n - a_n p_{n-1}.
\end{equation}
It is the coefficient $a_n$ herein which links with $\Delta_n$. Thus one has
\begin{equation}\label{an}
a_n^2 = \frac{\Delta_{n-1} \Delta_{n+1}}{\Delta_n^2}.
\end{equation}
Also relevant is the fact that the multiple integral
\begin{eqnarray}\label{Dw}
&& D_N(y_1,y_2)[w(x)] := \frac{1}{N!} \int_0^\infty dx_1 \cdots  \int_0^\infty dx_N \, \prod_{l=1}^N w(x_l) 
(y_1 - x_l) (y_2 - x_l)  \\ && \times \prod_{1 \le j < k \le N} (x_k - x_j)^2 \nonumber
\end{eqnarray}
can be expressed in terms of the polynomials $\{p_n\}$ according to the Christoffel-Darboux summation
(see e.g.~\cite{Fo02})
\begin{equation}\label{Dw1}
D_N(y_1,y_2)[w(x)] = \frac{\Delta_N}{\gamma_N \gamma_{N+1} }
\frac{p_{N+1}(y_1) p_N(y_2) - p_{N+1}(y_1) p_N(y_2)}{y_1 - y_2}.
\end{equation}

In random matrix theory the average (\ref{bc.3}) in the case $\mu=0$ is the generating function for the
probability that there are exactly $k$ eigenvalues in the interval $(0,t)$. In the case $\zeta = 0$
it gives the moment of the modulus of the characteristic polynomial $\prod_{l=1}^N(t - x_l)$. Moreover,
selecting the coefficient of $\zeta$ in \eqref{eq1:weight}, taking $\mu=2$ and substituting in
(\ref{Dw}) shows
\begin{eqnarray}\label{1.10}
&& D_N(s,s)[\theta(x-t) (x-t)^2 x^\alpha e^{-x}]  \\
&& \quad = \frac{1}{N!} \int_t^\infty dx_1 \cdots  \int_t^\infty dx_N \,
\prod_{l=1}^N (x_l - t)^2 (x_l - s)^2 x_l^\alpha e^{-x_l}
\prod_{j < k}^N (x_k - x_j)^2. \nonumber
\end{eqnarray}
After multiplying by
$$
e^{-t-s}(s-t)^2 (st)^\alpha
$$
this integral can be recognized as being proportional to the joint probability
density function of the first and second smallest eigenvalues, denoted $s$
and $t$ respectively, in the Laguerre unitary ensemble \cite{Forrester}. According to
(\ref{Dw1}) we can calculate (\ref{1.10}) in terms of quantities relating to the
polynomial system for the weight (\ref{eq1:weight}).

\subsection{Main results}
To present our main results requires some facts from the Okamoto $\tau$-function theory
of the fifth Painlev\'e equation \cite{Okamoto}.  The fifth Painlev\'e equation, 
\begin{eqnarray}\label{1.10a}
 && y'' = \Big ( \frac{1}{2y} + \frac{1}{y - 1} \Big ) (y')^2
- \frac{1}{t} y' + \frac{(y-1)^2}{t^2} \Big ( \alpha_1 y + \frac{\alpha_2}{y} \Big ) \\
&& \qquad + \frac{\alpha_3 y}{t} + \frac{\alpha_4 y (y+1)}{y - 1},  \nonumber
\end{eqnarray}
is one of the six non-linear differential equations identified by Painlev\'e and his students
as being distinct from classical equations and having the special property that all movable singularities 
are poles.
The Okamoto $\tau$-function theory relating to (\ref{1.10a})
is based on the Hamiltonian system with Hamiltonian specified by
\begin{eqnarray}\label{eq1:HV}
&& t H = q(q-1)^2 p^2 - \Big ( (v_2 - v_1) (q-1)^2 - 2(v_1 + v_2) q(q-1) + tq \Big ) p  \\
&& \qquad + (v_3 - v_1)(v_4 - v_1) (q - 1).  \nonumber 
\end{eqnarray}
Here $v_1,\dots,v_4$ are parameters constrained by
$$
v_1 + v_2 + v_3 + v_4 = 0.
$$
An essential feature of the theory is that eliminating $p$ in the Hamilton equations
$$
q' = \frac{\partial H}{\partial p}, \qquad p' = - \frac{\partial H}{\partial q}
$$
(the dashes denote differentiation with respect to $t$) gives that $q$ satisfies (\ref{1.10a})
with
$$
\alpha_1 = \frac{1}{2} (v_3 - v_4)^2, \quad
\alpha_2 = - \frac{1}{2} (v_2 - v_1)^2, \quad
\alpha_3 = 2v_1 + 2v_2 - 1, \quad
\alpha_4 = - \frac{1}{2}.
$$
The first of our main results may now be stated.

\begin{prop}\label{prop:maintheorem}
Let $v_1,\dots, v_4$ be such that
$$
v_3 - v_4 = - \mu, \quad v_3 - v_1 = n + \alpha, \quad v_3 - v_2 = n, \quad v_2 - v_1 = \alpha,
$$
and so
$$
\alpha_1 = \frac{\mu^2}{2}, \quad \alpha_2 = -\frac{\alpha^2}{2}, \quad
\alpha_3 = -(2n + \alpha + 1 + \mu), \quad \alpha_4 = - \frac{1}{2}
$$
In terms of the coefficients $\{a_n,b_n\}$ in the three term recurrence (\ref{eq0:3termrecurrence})
and the parameters $\alpha, \mu, t$ of the weight \eqref{eq1:weight}, let
\begin{subequations}\label{1.15}
\begin{eqnarray}
&& \theta_n = b_n - 2n - 1 - \alpha - \mu -t \\
&& \kappa_n = \left( n + \frac{\mu}{2} \right) t + a_n^2 - \sum_{i=0}^{n-1} b_i.
\end{eqnarray}
\end{subequations}
We have that the Hamilton equations are satisfied by
\begin{subequations}
\begin{eqnarray}
q &=& \frac{\theta_n+t}{\theta_n}\\
p &=& \frac{\theta_n(t + \alpha + \mu/2) - \kappa_n)}{ t(\theta_n + t)},
\end{eqnarray}
\end{subequations}
and so
\begin{subequations}
\begin{eqnarray}
\theta_n &=& \frac{t }{q-1}\\
\kappa_n &=& t(n +\alpha + \mu/2 - pq).
\end{eqnarray}
\end{subequations}
This characterization is made unique by the specification (\ref{3.19}) of the
small $t$ expansions of $\theta_n$ and $\kappa_n$.
\end{prop}

We can also express $\theta_n$ and $\kappa_n$ in terms of a solution of the fifth
Painlev\'e equation (\ref{1.10a}) with parameters different to those given in
Proposition \ref{prop:maintheorem}. For this we require the fact \cite{Okamoto} that \eqref{1.10a}
is formally unchanged upon the transformations
\begin{eqnarray*}
&& (\alpha,\beta,\gamma,\delta) \mapsto (-\beta,-\alpha,-\gamma,\delta) \quad y \mapsto \frac{1}{y}
\end{eqnarray*}

\commentout{
corresponding to the transformation 
\begin{eqnarray*}
&& (v_1,v_2,v_3,v_4) \mapsto (-v_3,-v_4,-v_1,-v_2), \qquad t \mapsto - t \\
&& q \mapsto \frac{1}{q}, \qquad pq \mapsto - pq + v_3 - v_1.
\end{eqnarray*}}

\begin{cor}
Suppose $v_1,\dots,v_4$ in (\ref{eq1:HV}) are such that
$$
v_2 - v_1 = - \mu, \quad v_3 - v_4 = \alpha, \quad v_2 - v_4 = n+\alpha+1, \quad v_2 - v_3 = n+1,
$$
and furthermore $t \mapsto - t$, so that $q$ satisfies (\ref{1.10a}) with
$$
\alpha_1 = \frac{\alpha^2}{2}, \quad \alpha_2 = - \frac{\mu^2}{2}, \quad
\alpha_3 = 2n + \alpha + 1 + \mu, \quad \alpha_4 = - \frac{1}{2}.
$$
(note that mapping $t \mapsto - t$ in  (\ref{1.10a}) is equivalent to mapping $\alpha_3 \mapsto - \alpha_3$).
In this case the Hamilton equations 
are satisfied by
\begin{subequations}
\begin{eqnarray}
q &=& \frac{\theta_n}{\theta_n+t}\\
p &=& \frac{(\theta_n + t)(\kappa_n - \mu/2 + \theta_n(1+2n+t+\alpha+\mu+\theta_n))}{ t \theta_n},
\end{eqnarray}
\end{subequations}
and so
\begin{subequations}
\begin{eqnarray}
\theta_n &=& \frac{t q}{1-q}\\
\kappa_n &=& tpq - \frac{t^2}{(1-q)^2} + \frac{t\left(t+ \frac{\mu}{2} - q\left(2n+\alpha+1+\frac{3\mu}{2}\right) \right)}{1-q}
\end{eqnarray}
\end{subequations}
\end{cor}

\section{Differential equations for orthogonal polynomial systems}

In the classical theory of orthogonal polynomials, the study of differential equations satisfied by orthogonal
polynomials has a long and distinguished history \cite{Shohat}. Under certain conditions \cite{Bonan}, given
a system of orthogonal polynomials, the derivatives of the polynomials may be expressed in terms of a linear
combination of at most two polynomials of the same system \cite{Laguerre, Bonan, IsmailChenLadder}. To 
describe these differential equations, we parameterize the coefficients of the polynomials according to
\begin{equation}\label{eq1:polycoeffs}
p_n(x) = \gamma_n x^n + \gamma_{n,1} x^{n-1} + \ldots + \gamma_{n,n}.
\end{equation}
In terms of this parameterization, it is clear from \eqref{eq0:3termrecurrence} that 
\begin{eqnarray*}
a_n &=& \frac{\gamma_{n-1}}{\gamma_n}\\
b_n &=& \frac{\gamma_{n,1}}{\gamma_n}-\frac{\gamma_{n+1,1}}{\gamma_{n+1}}.
\end{eqnarray*}

Our starting point is the condition that the logarithmic derivative of the weight is rational, and hence we write
\begin{equation}\label{eq2:logdrivweight}
\frac{\mathrm{d}}{\mathrm{d}x} \ln w(x) = \frac{2V}{W}
\end{equation}
where $W$ and $V$ are polynomials in $x$. This condition implies that the moments given by (\ref{eq0:moment}) 
satisfy some
recurrence relation. In particular, if we define the Stieltjes function by
\[
f(x) = \sum_{k = 0}^{\infty} \mu_k x^{-k-1} = \int_I \frac{w(s)}{x-s}ds
\]
then this recurrence is equivalent to $f$ satisfying the holonomic differential equation
\begin{equation}\label{eq2:holonomic}
W \frac{\mathrm{d}}{\mathrm{d}x} f = 2Vf + U
\end{equation}
where $U$ is some polynomial in $x$ of degree less than the degree of $V$. 

We now define the associated 
polynomials $\phi_{n-1}$ and associated functions $\epsilon_n$ by the equation
\begin{equation}\label{eq2:fpn}
fp_n = \phi_{n-1} + \epsilon_n.
\end{equation}
Explicitly
\begin{subequations}
\begin{eqnarray}
\label{eq1:epsint}\epsilon_n &=& \int_I \frac{p_n(s)}{x-s} w(s)\mathrm{d}s\\
\phi_{n-1} &=& \int_I \frac{p_n(s)-p_n(x)}{s-x} w(s)\mathrm{d}s
\end{eqnarray}
\end{subequations}
showing that $\phi_{n-1}$ is a polynomial of degree $n-1$
and $\epsilon_n$ is meromorphic at $x = \infty$. On this latter point, by orthogonality, 
\begin{equation}\label{eq1:epsasympty}
\epsilon_n \sim \gamma_n^{-1} x^{-n-1}
\end{equation}
as $x$ tends to $\infty$, while upon multiplying \eqref{eq0:3termrecurrence} by $f$ it is also
clear that the sequence of associated functions $\{\epsilon_n\}_{n=0}^{\infty}$ satisfies \eqref{eq0:3termrecurrence}. 
Using \eqref{eq0:3termrecurrence} and \eqref{eq1:polycoeffs}, we have the large $x$ expansions,
\begin{subequations}\label{eq2:largexexpansions}
\begin{eqnarray}
\label{eq2:largexpn}&& \quad p_n = \\
&& \nonumber \quad \gamma_n\left(x^n - x^{n-1}\sum_{i =0}^{n-1} b_i + x^{n-2}\left(\sum_{i =1}^{n-1} \sum_{j=0}^{i-1} b_i b_j -
\sum_{i=1}^{n-1}a_i^2\right) + O \left( x^{n-3} \right) \right) \\
\label{eq2:largexen}&& \quad \epsilon_n = \\
&&\hspace{-.1cm} \nonumber  \quad \gamma_n^{-1}\left(x^{-n-1} + x^{-n-2}\sum_{i =0}^n b_i + x^{-n-3}\left(\sum_{i =0}^{n} \sum_{j=0}^{i} b_i
b_j + \sum_{i=1}^{n+1}a_i^2\right) + O \left( x^{n-3} \right) \right).
\end{eqnarray}
\end{subequations}
Equating $(fp_n)p_{n-1}$ with $(fp_{n-1}) p_n$ shows $\phi_{n-1}p_{n-1} + \epsilon_np_{n-1} = \phi_{n-2}p_n + \epsilon_n$.
This and \eqref{eq2:largexexpansions} gives 
\begin{equation}\label{eq2:anrelation}
p_n\epsilon_{n-1} - p_{n-1}\epsilon_n = p_{n-1}\phi_{n-1} - p_n \phi_{n-2} = \frac{1}{a_n}.
\end{equation}
The above describes the notation and set formulae to be used in coming sections.

\subsection{Derivation via recurrence relations for moments and isomonodromy}\label{sec:recmoments}

In this section methods are outlined 
for obtaining differential equations satisfied by the 
orthogonal polynomial system based on the existence of a 
recurrence for the moments of the weight function, namely \eqref{eq2:holonomic}. 
This coupled with the theory of isomonodromic deformations allows us to construct the differential 
equations that govern evolution of the polynomials in both $x$ and $t$. 

\begin{thm}\label{thm:Magnus}
The orthogonal polynomials corresponding to a weight $w$ satisfy the differential equation 
\begin{equation}\label{eq2:difmagnus}
W(x)\frac{\mathrm{d}}{\mathrm{d}x}p_n(x) = (\Omega_n(x) - V(x))p_n(x) - a_n \Theta_n(x)p_{n-1}(x)
\end{equation}
where $\Omega_n$ and $\Theta_n$ are polynomials given by
\begin{subequations}
\begin{align}
& \Theta_n = W\left( \epsilon_n \frac{\mathrm{d}}{\mathrm{d}x} p_n - 
p_n\frac{\mathrm{d}}{\mathrm{d}x}\epsilon_n \right) + 2V \epsilon_np_n \label{eq2:Theta}\\
& \Omega_n = a_nW\left( \epsilon_{n-1} \frac{\mathrm{d}}{\mathrm{d}x}p_n - p_{n-1}\frac{\mathrm{d}}{\mathrm{d}x}\epsilon_n  \right) + a_nV(\epsilon_{n-1}p_n + \epsilon_np_{n-1})\label{eq2:Omega}.
\end{align}
\end{subequations}
\end{thm}

\begin{proof}
First we note from
\eqref{eq2:holonomic} and \eqref{eq2:fpn} that 
\[
W \frac{\mathrm{d}}{\mathrm{d}x}\left(\frac{\phi_{n-1}}{p_n}\right) -  \frac{2V\phi_{n-1}}{p_n} - U  = \frac{2V\epsilon_n}{p_n} - W\frac{\mathrm{d}}{\mathrm{d}x} \left(\frac{\epsilon_n}{p_n} \right). 
\]
Multiplying through by $p_n^2$ shows (\ref{eq2:Theta}) can be rewritten
\begin{equation}\label{2.10c}
\Theta_n = W\left( p_n \frac{\mathrm{d}}{\mathrm{d}x} \phi_{n-1} - \phi_{n-1} \frac{\mathrm{d}}{\mathrm{d}x}p_n \right ) - 2V\phi_{n-1}p_n - Up_n^2 .
\end{equation}
This tells us
that $\Theta_n$ is a polynomial, while (\ref{eq2:Theta}) bounds the degree. 
Explicitly, examining the $x \to \infty$ behaviour, namely \eqref{eq1:polycoeffs} and \eqref{eq1:epsasympty}, shows 
\begin{equation}\label{eq2:thetadegree}
\deg \Theta_n \leq \max( \deg W-2, \deg V - 2,0)
\end{equation}
as noted in \cite{Magnus1992}.

Using \eqref{eq2:anrelation} and (\ref{2.10c}) we find 
\begin{eqnarray*}
&& a_n (p_{n-1} \phi_{n-1} - p_n\phi_{n-2})\Theta_n = \\
&& \quad\quad  W\left( p_n \frac{\mathrm{d}}{\mathrm{d}x} \phi_{n-1} - \phi_{n-1} \frac{\mathrm{d}}{\mathrm{d}x}p_n \right ) - 2V\phi_{n-1}p_n - Up_n^2.
\end{eqnarray*}
By appropriately grouping terms divisible by $\phi_{n-1}$ and $p_n$ on opposite sides, we define the polynomial $\Omega_n$ to be the common factor according to 
\begin{eqnarray}
\label{eq2:diffOmega} p_n\phi_{n-1} \Omega_n &=& \phi_{n-1} \left( a_n\Theta_n p_{n-1} + W\frac{\mathrm{d}}{\mathrm{d}x}p_n  + V p_n \right) \\
&=& p_n \left( a_n \Theta_n \phi_{n-2} + W \frac{\mathrm{d}}{\mathrm{d}x}\phi_{n-1} - V\phi_{n-1} - Up_n \right).\nonumber
\end{eqnarray}
The first expression in \eqref{eq2:diffOmega} is equivalent to \eqref{eq2:difmagnus}
provided (\ref{eq2:Omega}) can be verified. 
For this purpose we use \eqref{eq2:Theta} in \eqref{eq2:diffOmega} to obtain
\[
\Omega_n = \frac{a_n W p_{n-1}\epsilon_n \frac{\mathrm{d}}{\mathrm{d}x} p_n}{p_n} - \frac{a_n W p_{n-1}p_n\frac{\mathrm{d}}{\mathrm{d}x}\epsilon_n}{p_n} + \frac{2V \epsilon_np_{n-1}p_n}{p_n} + \frac{W\frac{\mathrm{d}}{\mathrm{d}x}p_n}{p_n} + V.
\]
By rearranging \eqref{eq2:anrelation}, we find $a_n\epsilon_n p_{n-1} = a_n p_n \epsilon_{n-1}  -1$, which we use to remove occurrences of $p_{n-1}$, 
giving \eqref{eq2:Omega} as required. Examining the $x \to \infty$ behaviour, by using \eqref{eq1:polycoeffs} and \eqref{eq1:epsasympty} in \eqref{eq2:Omega}, shows 
\begin{equation}
\deg \Omega_n \leq \max( W-1, V)
\end{equation}
which again appears in \cite{Magnus1992}.
\end{proof}

The origin of this theorem can be traced back to the work of Laguerre \cite{Laguerre} and has since been revisted by contemporaries \cite{Ba90,BC90,Magnus1992}. The theorem provides a mechanical way of determining the differential equation satisfied by polynomials 
provided one knows the rational logarithmic derivative. One need only expand \eqref{eq2:Theta} and \eqref{eq2:Omega} to polynomial orders using \eqref{eq2:largexen} to 
produce a parametrization of the differential equation satisfied by the polynomials in terms of the $a_n$'s and $b_n$'s. A simple application of \eqref{eq0:3termrecurrence} gives us an expression for the derivative of $p_{n-1}$, which is but one column solution to 
a $2\times 2$ linear differential equation in $x$. The following corollary provides us with another solution. 

\begin{cor}
The function $\epsilon_n/w$ satisfies \eqref{eq2:difmagnus}.
\end{cor}
\begin{proof}
Consider the derivative of $fp_n$ in terms of $\epsilon_n$ and $\phi_n$. According to (\ref{eq2:fpn}) 
\[
W \frac{\mathrm{d}}{\mathrm{d}x} fp_n = W\frac{\mathrm{d}}{\mathrm{d}x}\phi_{n-1} + W\frac{\mathrm{d}}
{\mathrm{d}x}\epsilon_n.
\]
On the other hand, use of (\ref{eq2:difmagnus}) and (\ref{eq2:holonomic}) shows
\begin{eqnarray*}
W \frac{\mathrm{d}}{\mathrm{d}x} fp_n  &=& f \left( W \frac{\mathrm{d}}{\mathrm{d}x}p_n  \right)  + p_n W \frac{\mathrm{d}}{\mathrm{d}x} f \\
&=& (\Omega_n - V )fp_n -a_n \Theta_n fp_{n-1} + 2Vfp_n + Up_n\\
&=& (\Omega_n + V) \phi_{n-1} - a_n \Theta_n \phi_{n-1} + Up_n\\
&&\quad + (\Omega_n + V) \epsilon_n - a_n\Theta_n \epsilon_{n-1},
\end{eqnarray*}
where in obtaining the final equality (\ref{eq2:fpn}) has also been used.
By cancelling out the derivative of $\phi_{n-1}$ as calculated from the first expression for 
$\Omega_n$ from the previous proof, we deduce that the
derivative of $\epsilon_n$ is given by
\[
W \frac{\mathrm{d}}{\mathrm{d}x} \epsilon_n  = (\Omega_n + V) \epsilon_n - a_n \Theta_n \epsilon_{n-1}. 
\]
Hence
\begin{eqnarray*}
W \frac{\mathrm{d}}{\mathrm{d}x} \frac{\epsilon_n}{w} &=& \frac{w W \frac{\mathrm{d}}{\mathrm{d}x} \epsilon_n - \epsilon_n W \frac{\mathrm{d}}{\mathrm{d}x} w }{w^2} \\
&=& \frac{(\Omega_n + V) \epsilon_n - a_n \Theta_n \epsilon_{n-1} - 2V\epsilon_n}{w},
\end{eqnarray*}
where use has also been made of (\ref{eq2:logdrivweight}), as required.
\end{proof}

As a linear system, we have two linearly independent solutions. As mentioned above, both $\{p_n\}$ and $\{\epsilon_n/w\}$ satisfy \eqref{eq0:3termrecurrence} and \eqref{eq2:difmagnus}, telling us that $p_{n-1}$ and $\epsilon_{n-1}$ satisfy
\[
W \frac{\mathrm{d}}{\mathrm{d}x} y_n = a_n \Theta_{n-1}y_n + (\Omega_{n-1} -V - (x-b_n))y_{n-1}. 
\]
Hence, the matrix
\begin{equation}\label{2.13a}
Y_n = \begin{pmatrix}
p_n & \frac{\epsilon_n}{w}\\
p_{n-1} & \frac{\epsilon_{n-1}}{w}
\end{pmatrix}
\end{equation}
satisfies the matrix differential equation
\begin{equation}\label{eq2:linearx}
\frac{\mathrm{d}}{\mathrm{d}x}Y_n = \mathscr{A}_n Y_n
\end{equation}
where
\[
\mathscr{A}_n = \frac{1}{W}\begin{pmatrix} \Omega_n - V & -a_n \Theta_n \\ 
-a_n \Theta_{n-1} & \Omega_{n-1}-V - (x-b_{n-1})\Theta_{n-1} 
\end{pmatrix}.
\]
Because $\{ p_n \}$ and $\{\epsilon_n\}$ satisfy \eqref{eq0:3termrecurrence}, $Y_n$ also satisfies
\begin{equation}\label{eq2:linearn}
Y_{n+1} = M_n Y_n
\end{equation}
where
\[
M_n = \begin{pmatrix}
\frac{x-b_n}{a_{n+1}} & -\frac{a_n}{a_{n+1}} \\
1 & 0 
\end{pmatrix}.
\]

\begin{lem}
The polynomials $\Theta_n$ and $\Omega_n$ satisfy the recurrence relations
\begin{subequations}\label{eq2:Frueds}
\begin{eqnarray}
\label{eq2:Frued1}W + a_{n+1}^2\Theta_{n+1}- a_n^2 \Theta_{n-1} = (x-b_n)(\Omega_{n+1} - \Omega_n)\\
\label{eq2:Frued2}(x-b_{n-1})\Theta_{n-1} - (x-b_n)\Theta_n = \Omega_n - \Omega_{n+1}
\end{eqnarray}
\end{subequations}
\end{lem}

\begin{proof}
The equation \eqref{eq2:linearx} and \eqref{eq2:linearn} gives two ways of calculating $\frac{\mathrm{d}}{\mathrm{d}x}Y_{n+1}$. 
The consistency may be written as
\[
M_n \mathcal{A}_n - \mathcal{A}_{n+1} M_n + \frac{\mathrm{d}M_n}{\mathrm{d}x} = 0.
\]
This is an identity on the bottom two rows, however, in the first row, the consistency relation yeilds \eqref{eq2:Frued1} and \eqref{eq2:Frued2}
\end{proof}

This lemma can be found in Magnus \cite{Magnus1992}. Using \eqref{eq2:anrelation}, we have 
\[
\det Y_n  = \frac{1}{a_nw}.
\]
In general, for an equation of the form \eqref{eq2:linearx}, we see from (\ref{2.13a}) that 
\[
\frac{\mathrm{d}}{\mathrm{d}x} \det Y_n = \Tr \mathscr{A}_n \det Y_n
\]
and so
\[
\frac{\mathrm{d}}{\mathrm{d}x} \frac{1}{a_n w} = -\frac{2V}{w a_nW} = \Tr \mathscr{A}_n \det Y_n
\]
giving the additional relation
\begin{equation}\label{eq2:Frued3}
(x-b_n)\Theta_n = \Omega_{n+1} + \Omega_{n}. 
\end{equation}
This also implies \eqref{eq2:Frued2}. It further gives us a new parameterization of $\mathscr{A}_n$, given by
\[
\mathscr{A}_n = \frac{1}{W}\begin{pmatrix}
\Omega_n - V & -a_n\Theta_n \\
a_n \Theta_{n-1} & - \Omega_n - V
\end{pmatrix}
\]
as first derived in \cite{Magnus1992}. 

Another useful relation comes from the 
multiplication of \eqref{eq2:Frued1} and \eqref{eq2:Frued3}, which gives
\begin{equation}\label{2.17a}
W\Theta_n + a_{n+1}^2 \Theta_{n+1} \Theta_n - a_n^2 \Theta_n\Theta_{n-1}  = \Omega_{n+1}^2 - \Omega_n^2.
\end{equation}
Summing over $n$, given that $\Omega_0 = V$, shows
\begin{equation}\label{eq2:sumTheta}
\Omega_n^2 - a_n^2 \Theta_n \Theta_{n-1} = V^2 + W \sum_{i = 0}^{n-1} \Theta_i.
\end{equation}
 
\begin{figure}[!ht]
\begin{pspicture}(10,5)
\qdisk(3,3){.1}
\qdisk(6,3){.1}
\qdisk(9,3){.1}
\qdisk(4,1){.1}
\psbezier{->}(4,1)(1,5)(3,5)(4,1.1)
\psbezier{->}(4,1)(6,5)(10,5)(4.05,1.086)
\psbezier{->}(4,1)(10,5)(13,3)(4.086,1.05)
\rput(3,4){$\rho_1$}
\rput(6,4){$\rho_2$}
\rput(9,3.5){$\rho_3$}
\rput(3,3.2){$x_1$}
\rput(6,3.2){$x_2$}
\rput(8.7,3){$x_3$}
\end{pspicture}
\caption{\label{fig2:monodromy}The fundamental group of the complement of a set of three poles in $\mathbb{C}$.}
\end{figure}
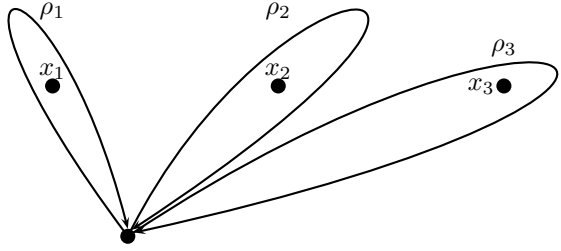

The roots of $W$ are now the poles of $\mathscr{A}_n$. If $\{x_j\}$ is the set of poles of $\mathscr{A}_i$, 
then we may write \eqref{eq2:linearx} as
\begin{equation}\label{eq2:monodromylinear}
\frac{\mathrm{d}}{\mathrm{d}x}Y_n = \mathscr{A}_nY_n = \left(\sum_i \frac{\mathscr{A}_{i,n}}{x- x_i} \right)Y_n
\end{equation}
where
\[
\mathscr{A}_{i,n} = \frac{1}{W'(x_i)}\begin{pmatrix}
\Omega_n(x_i) - V(x_i) & -a_n\Theta_n(x_i) \\
a_n \Theta_{n-1}(x_i) & - \Omega_n(x_i) - V(x_i)
\end{pmatrix}.
\]
This now places the differential equation into the context of isomodromic deformations. In general, any solution to \eqref{eq2:monodromylinear} is going to
be multivalued, with branch points at $\{x_i\}$(one possibly being $\infty$). Hence, by integrating around a path, say $\rho:[0,1]\to \mathbb{C}\setminus \{x_i\}$ where
$\rho(0)= \rho(1)$, the multi-valuedness can be expressed through the equation
\[
Y(\rho(0)) = Y(\rho(1)) \mathscr{M}_\rho
\]
where $\mathscr{M}_\rho$ is referred to as a monodromy matrix.
The set of monodromy matrices, $\{\mathscr{M}_\rho \}$, forms a representation of the fundemental group, $\pi_1(\mathbb{C}\setminus \{x_i\})$. This is depicted in figure \ref{fig2:monodromy}.

By construction, one solution of \eqref{eq2:monodromylinear} involves the polynomials, which are entire, and the associated functions. The goal of monodromy preserving deformations is to describe a family of linear problems of the form \eqref{eq2:monodromylinear} that share the same representation. A natural choice of deformation parameter turns out to be the poles of $\mathcal{A}_n$, giving rise to the classical result known as the Schlesinger equations, given by
\begin{subequations}\label{eq1:Schlesinger}
\begin{eqnarray}
\frac{\partial A_i}{\partial \alpha_j} &=& \frac{[A_i, A_j]}{\alpha_i - \alpha_j} \hspace{3cm} i \neq j\\
\frac{\partial A_i}{\partial \alpha_i} &=& -\sum_{j \neq i} \frac{[A_i, A_j]}{\alpha_i - \alpha_j}.
\end{eqnarray}
\end{subequations}
We shall assume that just one of the poles depends on a variable, $t$, which shall become the deformation parameter. 
In the case that we have just one parameter that needs to be deformed, we have that there is a matrix 
\begin{equation}\label{eq2:lineart}
\mathscr{B}_n(x,t) = \frac{\partial Y_n}{\partial t} Y_n^{-1}
\end{equation}
The form of this matrix, as implied by \eqref{eq1:Schlesinger}, is given by \cite{Chud}
\begin{equation}\label{2.21'}
\mathscr{B}_n = \mathscr{B}_{\infty,n} - \sum_i \frac{\mathscr{A}_{i,n}}{x - x_i} \frac{\partial x_i}{\partial t}.
\end{equation}
By examining the large $x$ behaviour of $p_n$ and $\frac{\partial p_i}{\partial t}$, we deduce that $\mathscr{B}_{\infty,n}$ is given by
\begin{equation}\label{2.21a}
\mathscr{B}_{\infty,n} = \begin{pmatrix}
\frac{1}{\gamma_n}\frac{\partial \gamma_n}{\partial t} & 0 \\
0 &- \frac{1}{\gamma_{n-1}}\frac{\partial \gamma_{n-1}}{\partial t}
\end{pmatrix}.
\end{equation}
This gives us two linear differential equations for the one system, which must be consistent. Hence we have the condition
\[
\frac{\partial}{\partial t}\frac{\partial}{\partial x} Y_n = \frac{\partial}{\partial x}\frac{\partial}{\partial t} Y_n,
\]
which is equivalent to
\begin{equation}\label{eq1:Comp}
\mathscr{A}_n \mathscr{B}_n - \mathscr{B}_n \mathscr{A}_n + 
\frac{\partial}{\partial t} \mathscr{A}_n - \frac{\partial}{\partial x}
\mathscr{B}_n = 0.
\end{equation}
This completely determines the differential equation for the orthogonal polynomials and associated functions in $t$.

\subsection{Ladder Operators}\label{sec:ladder}

An alternative approach, developed by Chen and collaborators
\cite{IsmailChenLadder, Chen:discontinuous, BasorChen}, is that of the ladder operators. 
We shall assume that the weight satisfies the same logarithmic differential equation 
(\ref{eq2:logdrivweight}), and hence that the corresponding moments satisfy (\ref{eq2:holonomic}).
However, this approach typically concerns the monic versions of the orthgonal polynomials, given by
\[
P_n = \frac{1}{\gamma_n}p_n.
\]
In order to make comparisons to the previous section, we will deal primarily with $p_n$ rather than $P_n$.
Now $p_n$, being a polynomial of degree $n$, when differentiated can be expressed as a linear
combination of $\{p_j\}_{j=0,\dots,n-1}$,
\begin{equation}\label{2.22a}
\frac{\mathrm{d}}{\mathrm{d}x}p_n = \sum_{i =0}^{n-1} \alpha_{n,i} p_i.
\end{equation}
We may reduce this, via the use of \eqref{eq0:3termrecurrence}, to a differential equation specified by the following
theorem of Bonan and Clark \cite{BC90} and Bauldry \cite{Ba90}.

\begin{thm}\label{thm:ladder}
The orthogonal polynomial system $\{ p_n \}$, defined by the weight $w$, satisfies 
\begin{equation}\label{eq2:ladder}
\frac{\mathrm{d}}{\mathrm{d}x}p_n = - B_n p_n + a_n A_n p_{n-1}\\
\end{equation}
where 
\begin{subequations}\label{eq2:laddercoefficients}
\begin{eqnarray}
A_n &=& \int_I \left( \frac{p_n(y)^2}{x-y} \right) \left( \frac{2V(x)}{W(x)} - \frac{2V(y)}{W(y)}\right) w(y) \\
\label{eq2:ladderBdef}B_n &=& a_n \int_I  \left( \frac{p_n(y)p_{n-1}(y)}{x-y}\right) \left( \frac{2V(x)}{W(x)} - \frac{2V(y)}{W(y)}\right) w(y) \mathrm{d}y.
\end{eqnarray}
\end{subequations}
\end{thm}

\begin{proof}
Beginning with (\ref{2.22a}),
we may use orthogonality and integration by parts to find 
\[
\alpha_{n,i} = \int_I p_n'(y) p_i(y) w(y) \mathrm{d}y = - \int_I p_n(y) \left( p_i'(y) - \frac{2p_i(y)V}{W}\right) w(y) \mathrm{d}y.
\]
However, since $p_i'$ is a polynomial whose degree is less than $n$, this term must be destroyed by orthogonality, leaving
\[
\alpha_{n,i} = -\int_I \frac{2V}{W}p_n(y) p_k(y) w(y) \mathrm{d}y.
\]
The derivative of $p_n$ can therefore be written 
\[
p_n'(x) = - \int_I  \sum_{i=0}^{n-1} p_i(x)p_i(y) \frac{2V(y)}{W(y)}w(y)\mathrm{d}y.
\]
However, if we replace $2V(y)/W(y)$ with $2V(x)/W(x)$, this would vanish by orthogonality, hence we may add it to obtain
\[
p_n'(x) = \int_I p_n(y) \left( \sum_{i = 0}^{n-1} p_i(x)p_i(y) \right)\left( \frac{2V(x)}{W(x)} - \frac{2V(y)}{W(y)}\right) w(y) \mathrm{d}y.
\]
By exploiting the Christoffel-Darboux summation (\ref{Dw1})
and pulling out the polynomials in $x$, we arrive at
\begin{eqnarray*}
&& p_n'(x) =a_n p_{n-1}(x) \int_I  \left( \frac{p_n(y)^2}{x-y}\right) \left( \frac{2V(x)}{W(x)} - \frac{2V(y)}{W(y)}\right) w(y) \mathrm{d}y \\ 
&& \quad \quad - a_n p_n(x) \int_I \left( \frac{p_n(y)p_{n-1}(y)}{x-y} \right) \left( \frac{2V(x)}{W(x)} - \frac{2V(y)}{W(y)}\right) w(y) \mathrm{d}y ,
\end{eqnarray*}
which is \eqref{eq2:ladder}.
\end{proof}

This allows us to define the ladder operator,
\[
L_{n,1} := \left( \frac{\mathrm{d}}{\mathrm{d}x} + B_n \right)
\]
which has the effect
\[
L_{n,1} p_n = A_n p_{n-1}.
\]

\begin{lem}
The terms $A_n$ and $B_n$ satisfy the recurrence relations
\begin{subequations}
\begin{eqnarray}
\label{eq2:ladderrecurrence1}B_{n+1} + B_n &=& (x-b_n)A_n - \frac{2V}{W}\\
\label{eq2:ladderrecurrence2}(B_{n+1} - B_n)(x-b_n) &=& a_{n+1}^2 A_{n+1} -a_n^2 A_{n-1} + 1
\end{eqnarray}
\end{subequations}
\end{lem}

\begin{proof}
Using \eqref{eq2:ladderBdef},
\begin{eqnarray*}
B_n + B_{n+1}\hspace{-.3cm}  &=& \hspace{-.3cm} \int_I  \left( \frac{p_n(y)(a_n p_{n-1}(y) + a_{n+1}p_{n+1}(y))}{x-y}\right) 
\left( \frac{2V(x)}{W(x)} - \frac{2V(y)}{W(y)}\right) w(y) \mathrm{d}y\\
&=& \int_I  \left( \frac{p_n(y)^2(y-b_n)}{x-y}\right)  \left( \frac{2V(x)}{W(x)} - \frac{2V(y)}{W(y)}\right) w(y) \mathrm{d}y \\
&=& (x-b_n)A_n  + \int_I p_n(y)^2 \left( \frac{2V(x)}{W(x)} - \frac{2V(y)}{W(y)}\right) w(y) \mathrm{d}y\\
&=& (x-b_n)A_n  - \frac{2V}{W}.
\end{eqnarray*}
By consistency of \eqref{eq2:ladder} with \eqref{eq0:3termrecurrence} we obtain the second required expression.
\end{proof}

By multiplying and rearranging \eqref{eq2:ladderrecurrence1} and \eqref{eq2:ladderrecurrence2} we obtain 
\begin{equation}\label{2.17b}
B_{n+1}^2 - B_n^2 - \frac{2V}{W}(B_{n+1}- B_n) = a_{n+1}^2A_{n+1}A_n - a_n^2 A_{n-1}A_n + A_n.
\end{equation}
Hence, by summing over $n$ and appropriately evaluate initial conditions, we obtain
\begin{equation}\label{2.17c}
B_n^2 - \frac{2V}{W}B_n - a_n^2 A_n A_{n-1} = -\sum_{i=0}^{n-1} A_i.
\end{equation}
Note the structual correspondence of (\ref{2.17b}) and (\ref{2.17c}) with
(\ref{2.17a}) and (\ref{eq2:sumTheta}) respectively.

\section{Derivations of $\mathrm{P}_{V}$}

We now turn to the polynomials specified by \eqref{eq0:bilinear} with the weight specified by \eqref{eq1:weight}. 
The above formulae for the derivatives in $x$ should be considered as partial derivatives. 
The formula (\ref{eq0:moment}) for the moments is a hypergeometric integral, which may be evaluated to give
\begin{eqnarray}\label{eq2:momentsf}
 \mu_k &=& (\zeta - 1)\Gamma(1+k+\alpha+\mu){}_1 F_1 \left( \begin{array}{ c |} -\mu \\ -k-\alpha-\mu \end{array} \hspace{.2cm} -t\right)+\\
&&  \left( (-1)^\mu + \frac{(\zeta-1)\sin(\pi (\alpha))}{\sin(\pi(\alpha +\mu))}\right)\frac{\Gamma(k+1)\Gamma(\mu+1)!}{\Gamma(2+k+\alpha+\mu)} \times \nonumber \\
&&  t^{1+k+\alpha + \mu}{}_1F_1\left( \begin{array}{ c |} 1+k+\alpha \\ 2+k+\alpha+\mu \end{array} \hspace{.2cm} -t\right)\nonumber \\
&=& C_1(\zeta,\mu,k,\alpha) {}_1 F_1 \left( \begin{array}{ c |} -\mu \\ -k-\alpha-\mu \end{array} \hspace{.2cm} -t\right) \nonumber \\
&& + C_2(\zeta,\mu,k,\alpha) t^{1+k+\alpha + \mu} {}_1F_1\left( \begin{array}{ c |} 1+k+\alpha \\ 2+k+\alpha+\mu \end{array} \hspace{.2cm} -t\right) \nonumber 
\end{eqnarray}
where
\begin{eqnarray*}
C_1(\zeta,\mu,k,\alpha) &=& (1-\zeta)\Gamma(1+k+\alpha+\mu)\\
C_2(\zeta,\mu,k,\alpha) &=&  \left( (-1)^\mu + \frac{(\zeta-1)\sin(\pi \alpha)}{\sin(\pi(\alpha +\mu))}\right)\frac{\Gamma(k+1)\Gamma(\mu+1)}{\Gamma(2+k+\alpha+\mu)}
\end{eqnarray*}
and ${}_1F_1$ is the confluent hypergeometric function. We seek the corresponding differential equations 
satisfied by the orthogonal polynomial system, as implied by the theory of
\S \ref{sec:recmoments} and \S \ref{sec:ladder}

To derive the differential equation satisfied by the orthogonal polynomials from the theory of \S \ref{sec:recmoments},  we remark that since the factor of $(1 + \zeta\theta(x-t))$ plays the role of a multiplicative constant almost everywhere, the logarithmic derivative
of $w$ coincides with the logarithmic derivative of $\newblock{(x - t)^\mu x^\alpha e^{-x}}$ almost everywhere. Hence
we write
\[
x(x-t)\frac{\partial_x w}{w} \cong \left( -x^2 + (\alpha + \mu +t)x - \mu t\right)  
\]
where $\cong$ is to be interpreted as equals almost everywhere, and so independent of $\zeta$
\begin{eqnarray*}
W &=& x(x-t) \\
2V &=& -x^2 + (\alpha + \mu -t)x + \mu t.
\end{eqnarray*}
Recall that the form of the logarithmic derivative is the essential ingredient 
in both theorem \ref{thm:Magnus} and theorem \ref{thm:ladder}. 

\subsection{Recurrence of moments approach}

Now that $W$ and $V$ have been defined, determining the differential equation satisfied for this particular family of orthogonal polynomials is simply a matter of applying theorem \ref{thm:Magnus}.

\begin{cor}
The polynomials $p_n$ corresponding to the weight  \eqref{eq1:weight}
satisfy the differential equation
\begin{equation}\label{eq3:linearmag}
\partial_x \begin{pmatrix} p_n \\ p_{n-1} \end{pmatrix} =  \left\{ \mathscr{A}_\infty + \frac{\mathscr{A}_0}{x} + \frac{\mathscr{A}_t}{x-t}\right\} \begin{pmatrix} p_n \\p_{n-1} \end{pmatrix}
\end{equation}
where
\begin{eqnarray*}
\mathscr{A}_0 &=& \frac{1}{t} \begin{pmatrix} \kappa_n - \frac{\mu t}{2} & -a_n \theta_n \\
a_n \theta_{n-1} & - \kappa_n - \frac{\mu t}{2} \end{pmatrix}\\
\mathscr{A}_t &=& \frac{1}{t} \begin{pmatrix} \left(n+\frac{\mu}{2}\right)t - \kappa_n  & a_n (\theta_n +t)\\
- a_n(\theta_{n-1}+t) & \kappa_n - \left(n+\alpha+\frac{\mu}{2}\right)t \end{pmatrix}\\
\mathscr{A}_\infty &=& \begin{pmatrix} 0 & 0 \\ 0 & 1\end{pmatrix} 
\end{eqnarray*}
and
\begin{eqnarray*}
\theta_n &=& b_n - 2n - 1 - \alpha - \mu - t\\
\kappa_n &=& \left(n+\frac{\mu}{2}\right)t + a_n^2 + \frac{\gamma_{n,1}}{\gamma_n}.
\end{eqnarray*}
\end{cor}

\begin{proof}
By way of application of \eqref{eq2:Omega} and \eqref{eq2:Theta} using \eqref{eq2:largexexpansions} one obtains 
for the explicit form of $\Omega_n$ and $\Theta_n$, 
\begin{eqnarray*}
\Omega_n &=& -\frac{x^2}{2} + x \left(\frac{2n +\alpha + \mu + t}{2}\right) - \frac{2n+\mu t}{2} - a_n^2 - \frac{\gamma_{n,1}}{\gamma_n} \\
\Theta_n &=& -x + 2n+1+\alpha+\mu + t - b_n,
\end{eqnarray*}
which we subsequently decompose into the form seen above.
\end{proof}

We also require these equations to be written in terms of the $\kappa_n$ and $\theta_n$ alone. For this, we 
note that  recurrence relations for $\theta_n$ and $\kappa_n$
are implied by \eqref{eq2:Frued1} and \eqref{eq2:Frued3}. 

\begin{cor}\label{c3.4}
The associated functions, $\theta_n$ and $\kappa_n$, satisfy the recurrences
\begin{subequations}\label{eq2:kappathetarecurrences}
\begin{eqnarray}
\label{eq3:kq1}\kappa_{n+1} + \kappa_n &=& -\theta_n(\theta_n + t + 2n+\alpha+1 + \mu)\\
\label{eq3:kq2}\frac{\theta_n}{\theta_n + t} \frac{\theta_{n-1}}{\theta_{n-1} + t} &=& \frac{\kappa_n^2 - \frac{\mu^2 t^2}{4}}{\left(\kappa_n - 
\left(n+\alpha + \frac{\mu}{2}\right)t \right)\left(\kappa_n - \left(n+ \frac{\mu}{2}\right)t\right)}.
\end{eqnarray}
\end{subequations}
\end{cor}

\begin{proof}
The relation \eqref{eq2:Frued2} is equivalent to \eqref{eq3:kq1} when one uses the definitions of $\Theta_n$ and $\Omega_n$ in terms of $\theta_n$ and $\kappa_n$. Evaluating \eqref{eq2:sumTheta} at $x=0$ and $x=t$ shows
\begin{eqnarray}
\label{eq2:thetan-1intermsof}a_n^2\theta_n \theta_{n-1}&=& \kappa_n^2 - \frac{\mu^2 t^2}{4} \\
a_n^2(t+\theta_n)(t+ \theta_{n-1}) &=& \left(\kappa_n - \frac{(2n+2\alpha+\mu)t}{2} \right) \left(\kappa_n - \frac{(2n+\mu)t}{2} \right)
\end{eqnarray}
respectively. The ratio of these identities is \eqref{eq3:kq2}.
\end{proof}

The relation \eqref{eq3:kq2} may be used to eliminate the occurence of $\theta_{n-1}$ in \eqref{eq3:linearmag}. Note that \eqref{eq3:linearmag} is the form of \eqref{eq2:monodromylinear}, where only one of the poles depends on $t$. Hence, the evolution in $t$ is governed by \eqref{eq2:lineart}. In this regard, the derivations of the time derivatives for $\theta_n$ and $\kappa_n$ and the methods of Forrester and Witte \cite{Forrester} contrast with the methods of Basor and Chen \cite{BasorChen}. Once the derivatives in $x$ are found, one may apply the theory of isomondromic deformations \cite{JimboMiwaI} to obtain appropriate derivatives in $t$. This approach, as seen in \cite{Forrester} extends the evolution of the orthogonal polynomials in the $t$ direction via the following result.

\begin{cor}
In addition to (\ref{eq3:linearmag}), the othogonal polynomials satisfy
\begin{equation}\label{eq2:linearprobtf}
\partial_t \begin{pmatrix} p_n \\ p_{n-1} \end{pmatrix} = \left\{ \mathscr{B} - \frac{\mathscr{A}_t}{x-t}\right\}\begin{pmatrix} p_n \\ p_{n-1} \end{pmatrix}
\end{equation}
where
\[
\mathscr{B} = \frac{1}{2t}\begin{pmatrix} \theta_n + t & 0 \\ 0 & -\theta_{n-1} - t \end{pmatrix}.
\]
\end{cor}

\begin{proof}
This almost directly follows from the corollary of the Schlesinger equations, (\ref{c3.4}) and
(\ref{2.21'}). We use the fact that in the context of orthogonal polynomials $\mathscr{B}$ has the 
explicit form (\ref{2.21a}).
By equating the residues of the left and right hand side of 
the compatibility relation \eqref{eq1:Comp} at $x= \infty$, the diagonal entries reveal
\begin{subequations}\label{eq2:thetagammarelation}
\begin{eqnarray}
\frac{2\partial_t \gamma_n}{\gamma_n} &=& 1+ \frac{\theta_n}{t} \\
\frac{2 \partial_t \gamma_{n-1}}{\gamma_{n-1}}&=& 1+ \frac{\theta_{n-1}}{t}
\end{eqnarray}
\end{subequations}
while the off diagonal entries are $0$. This gives the required form for $\mathscr{B}$ above.
\end{proof}

One may easily calculate the derivatives in $t$ of $\theta_n$ and $\kappa_n$ via the compatibility of 
\eqref{eq2:linearx} and \eqref{eq2:lineart}, 
\begin{equation}\label{eq2:compatibility}
\partial_t  \mathscr{A} - \partial_x  \mathscr{B} +  \mathscr{A}  \mathscr{B} - 
 \mathscr B  \mathscr A = 0,
\end{equation}
to define the evolution of $\theta_n$ and $\kappa_n$. Using the above recursion relations allows one to express
the derivatives of $\theta_n$ and $\kappa_n$ in terms of themselves.
Alternatively, using the general framework of \cite{Magnus1992}, the derivatives of $a_n$ and $b_n$ are expressible 
in terms of the  functions $\Theta_n$ and $\Omega_n$ evaluated 
at the movable finite singular points of \eqref{eq2:linearx} 
via the expression
\begin{subequations}\label{eq2:magnustricktderivatives}
\begin{eqnarray}
\frac{\mathrm{d}}{\mathrm{d}t}\ln a_n &=& \frac{1}{2} \sum_{r = 1}^m \frac{\Theta_n(x_r) - \Theta_{n-1}(x_r)}{W'(x_r)} \frac{\mathrm{d}}{\mathrm{d}t}x_r\\
\frac{\mathrm{d}}{\mathrm{d}t}\ln b_n &=& \sum_{r = 1}^m \frac{\Omega_{n+1}(x_r) - \Omega_{n-1}(x_r)}{W'(x_r)} \frac{\mathrm{d}}{\mathrm{d}t}x_r,
\end{eqnarray}
\end{subequations}
where in the case of \eqref{eq2:linearx}, $m=1$, and the only point is $x_1 = t$. This leads to the equations
\begin{subequations}\label{eq2:abdifferential}
\begin{eqnarray}
\frac{2t}{a_n} \frac{\mathrm{d}a_n}{\mathrm{d}t} &=& 2 + b_{n-1} - b_n \\
t  \frac{\mathrm{d}b_n}{\mathrm{d}t} &=& a_n^2 - a_{n+1}^2 + b_n.
\end{eqnarray}
\end{subequations}
We know \eqref{eq2:magnustricktderivatives} is equivalent to \eqref{eq2:abdifferential}. Using \eqref{eq2:compatibility} 
and \eqref{eq2:abdifferential} in conjunction with \eqref{eq2:kappathetarecurrences} to eliminate occurences of $a_n^2$ 
and $\theta_{n-1}$ gives a differential system for $\{\theta_n,\kappa_n\}$. 

\begin{cor}\label{cor2:tderivatives}
The associated functions, $\theta_n$ and $\kappa_n$ satisfy the coupled differential equations 
in $t$
\begin{subequations}\label{eq2:forrestersthetarel}
\begin{eqnarray}
\label{eq2:forrestersthetarela}&t \frac{\partial}{\partial t} \theta_n = 2\kappa_n + (2n + \alpha+1+\mu+t+ \theta_n)\theta_n\\
\label{eq2:forrestersthetarelb}&\hspace{.2cm}t \frac{\partial \kappa_n}{\partial t}  = \left( \frac{1}{\theta_n + t} + \frac{1}{\theta_n}\right) \kappa_n^2 + \left(2n + \alpha + \mu+ 1 - (2n+ \alpha+ \mu)\frac{t}{\theta_n + t} \right) \kappa_n \\ 
& - \left(n^2 + \left(n+\frac{\mu}{2}\right)(\alpha + \mu)\right)t - \frac{\mu^2 t^2}{4\theta_n} + \left(n+\frac{\mu}{2}\right)\left(n+ \alpha+ \frac{\mu}{2}\right) \frac{t^2}{\theta_n + t}\nonumber.
\end{eqnarray}
\end{subequations}
\end{cor}

\begin{proof}
This simply follows from the evaluation of \eqref{eq1:Comp}. 
The first relation follows from \eqref{eq2:thetagammarelation}, namely
\[
\frac{2t a_n'}{a_n} = \theta_{n-1} - \theta_n.
\]
The two other relations that arise are
\begin{eqnarray*}
t\theta_n' &=& 2\kappa_n + (2n + \alpha+1+\mu+t+ \theta_n)\theta_n\\
t\kappa_n' &=&  \kappa_n - a_n^2 (\theta_n - \theta_{n-1}).
\end{eqnarray*}
The first of these is \eqref{eq2:forrestersthetarela}. By using \eqref{eq2:thetan-1intermsof} to eliminate $a_n^2$ and \eqref{eq3:kq2} to eliminate $\theta_{n-1}$, one obtains \eqref{eq2:forrestersthetarelb}
\end{proof}

Now that one has the derivatives of $\kappa_n$ and $\theta_n$, the remaining task is to find the
transformation which allows them to be identified as the Hamilton equations for a Painlev\'e V
system. But before doing this, we want to show how differential equations equivalent to the
coupled system (\ref{eq2:forrestersthetarel}) can be derived from the formalism of \S 2.2. 

\subsection{Ladder operator approach}

We want to specialize
\eqref{eq2:ladder} to the weight \eqref{eq1:weight}. 

\begin{prop}
With the weight (\ref{eq1:weight}) the coefficients $A_n$ and $B_n$ in \eqref{eq2:ladder} are given by
\begin{eqnarray*}
A_n &=& \frac{R_n}{x-t} + \frac{1-R_n}{x}\\
B_n &=& \frac{r_n}{x-t} - \frac{n + r_n}{x}
\end{eqnarray*}
where for $\alpha \geq 1$
\begin{eqnarray*}
R_n &=& \alpha \int_0^{\infty} \frac{w(y)p_n(y)^2}{(t-y)}\mathrm{d}y \\
r_n &=& \alpha \int_0^{\infty} \frac{w(y)p_n(y)p_{n-1}}{(t-y)}\mathrm{d}y
\end{eqnarray*}
\end{prop}

\begin{proof}
We note that
\[
\frac{2V(x)}{W(x)} - \frac{2V(y)}{W(y)} = \left(-\frac{\alpha}{(t-x)(t-y)} - \frac{\mu}{xy} \right)(x-y),
\]
and hence the integrals that define $A_n$ and $B_n$ simplify to 
\begin{eqnarray*}
A_n &=& \int_0^{\infty} p_n(y)^2 \left(-\frac{\alpha}{(t-x)(t-y)} - \frac{\mu}{xy} \right) w(y) \mathrm{d}y\\
B_n &=& \int_0^{\infty} p_n(y)p_{n-1}(y) \left(-\frac{\alpha}{(t-x)(t-y)} - \frac{\mu}{xy} \right) w(y) \mathrm{d}y,
\end{eqnarray*}
or equivalently
\begin{eqnarray*}
A_n &=& \frac{\alpha}{x-t} \int_0^{\infty} (1 - \zeta \theta(y - t) ) (y - t)^{\alpha-1} y^\mu e^{-y} p_n(y)^2 \mathrm{d}y\\
&& \hspace{1cm} - \frac{\mu}{x} \int_0^{\infty} (1 - \zeta \theta(y - t) ) (y - t)^\alpha y^{\mu-1} e^{-y} p_n(y)^2 \mathrm{d}y \\
B_n &=& \frac{\alpha}{x-t} \int_0^{\infty} (1 - \zeta \theta(y - t) ) (y - t)^{\alpha-1} y^\mu e^{-y} p_n(y)p_{n-1}(y) \mathrm{d}y \\ 
&& \hspace{1cm} - \frac{1}{x}\int_0^{\infty} (1 - \zeta \theta(y - t) ) (y - t)^\alpha y^{\mu-1} e^{-x} p_n(y)p_{n-1}(y) \mathrm{d}y. 
\end{eqnarray*}
Now define
\begin{eqnarray*}
R_n &=& \alpha \int_0^{\infty} (1 - \zeta \theta(y - t) ) (y - t)^{\alpha-1} y^\mu e^{-y} p_n(y)^2 \mathrm{d}y \\
&=& \alpha \int_0^{\infty} \frac{w(y)p_n(y)^2}{(t-y)}\mathrm{d}y \\
r_n &=& \alpha \int_0^{\infty} (1 - \zeta \theta(y - t) ) (y - t)^{\alpha-1} y^\mu e^{-y} p_n(y)p_{n-1}(y) \mathrm{d}y\\
&=& \alpha \int_0^{\infty} \frac{w(y)p_n(y)p_{n-1}}{(t-y)}\mathrm{d}y.
\end{eqnarray*}
We apply integration by parts, orthogonality and the known value of $w(x,t)$ at $0,t$ and $\infty$, 
to express the second part of the integrals in $A_n$ and $B_n$ in terms of $R_n$ and $r_n$ respectively, giving 
\begin{eqnarray*}
A_n &=& \frac{R_n}{x-t} - \frac{1}{x} \int_0^{\infty} \left(\frac{\alpha w(y,t)p_n^2(y)}{t-y} - w(y,t)p_n^2(y) \right)\mathrm{d}y \\
B_n &=& \frac{r_n}{x-t} - \frac{\mu}{x}\int_0^{\infty} \left(\frac{\alpha w(y,t)p_n(y)p_{n-1}}{t-y} + w(y,t) p_{n-1}(y) \frac{\partial p_n}{\partial y} \right) \mathrm{d}y .
\end{eqnarray*}
Using orthogonality and the expression
\[
\frac{\partial p_n(y)}{\partial y} = \frac{n}{a_n} p_{n-1}(y) + \textrm{ lower order terms } 
\]
gives the stated formulas.
\end{proof}

Now that the form of $A_n$ and $B_n$ are known, the differential equation satisfied by the polynomials can
be written 
\[
\frac{\mathrm{d}}{\mathrm{d}x} \begin{pmatrix} p_n \\ p_{n-1} \end{pmatrix} = \begin{pmatrix} \frac{n+r}{x} - \frac{r_n}{x-t} & a_n\left( \frac{1- R_n}{x} + \frac{R_n}{x-t}\right) \\
-a_n\left( \frac{1-R_{n-1}}{x} + \frac{R_{n-1}}{x-t}\right) & \frac{r_n}{x-t} - \frac{n+r_n}{x} - \frac{2V}{W} \end{pmatrix}\begin{pmatrix} p_n \\ p_{n-1} \end{pmatrix}
\]
where the second row is a consequence of \eqref{eq0:3termrecurrence} and \eqref{eq2:ladderrecurrence1}.
\begin{lem}
The functions $R_n$ and $r_n$ satisfy the recurrences
\begin{subequations}
\begin{eqnarray}
\label{eq3:laddrec1}r_{n+1} + r_n - \alpha &=& R_n(\mu + \alpha + 2n + 1 +tR_n - t)\\
\label{eq3:ladderc2}\frac{R_nR_{n-1}}{(R_n-1)(R_{n-1} - 1)} &=& \frac{r_n(r_n - \alpha)}{(r_n + n)(r_n + n+ \mu)}
\end{eqnarray}
\end{subequations}
\end{lem}

\begin{proof}
The residue in $x$ of \eqref{eq2:ladderrecurrence1} at $t$ and $\infty$ using these definitions for $A_n$ and $B_n$ in terms of $R_n$ and $r_n$ shows
\begin{eqnarray}
\label{eq3:bnintermsofRn}b_n &=& 2n+ 1 + \alpha + \mu + tR_n\\
r_{n+1} + r_n -\alpha &=& R_n (t-b_n)\nonumber
\end{eqnarray}
which gives \eqref{eq3:laddrec1}. The evaluation of the result of multiplying \eqref{eq2:ladderrecurrence2} by $x^2(x-t)^2$ at $0$ and $t$ reveals 
\begin{eqnarray}
\label{eq3:rnintermsofRnRnm1}r_n(r_n-\alpha) &=& a_n^2R_{n-1}R_n\\
(n + r_n)(n+\mu+r_n) &=& a_n^2(R_n-1)(R_{n-1}-1) \nonumber
\end{eqnarray}
giving \eqref{eq3:ladderc2} in an analogous manner to the previous section.
\end{proof}

In addition to \eqref{eq3:bnintermsofRn}, there is a further
relation obtained by eliminating $R_{n-1}$ from \eqref{eq3:rnintermsofRnRnm1} by using \eqref{eq3:ladderc2}, giving 
\begin{equation}\label{eq3:anintermsofRnandrn}
a_n^2 = \frac{(r_n-\alpha)r_n}{R_n} - \frac{(n+r_n)(n+\mu+r_n)}{R_n-1}.
\end{equation}

\begin{lem}
The recursion coefficients of \eqref{eq0:3termrecurrence} satisfy the differential equations
\begin{subequations}
\begin{eqnarray} 
\label{eq3:Ldan}\frac{2}{a_n}\frac{\mathrm{d}a_n}{\mathrm{d}t} &=& R_{n-1} - R_n \\
\label{eq3:Ldbn}\frac{\mathrm{d}b_n}{\mathrm{d}t} &=& r_n - r_{n+1}
\end{eqnarray}
\end{subequations}
\end{lem}

\begin{proof}
We take the derivative of \eqref{eq0:orthonormality} in the case $i= j = n$ to see that
\[
0 = \frac{\mathrm{d}}{\mathrm{d}t}\int_{0}^{\infty} p_n(y)^2 w(y) \mathrm{d}y 
\]
with respect to $t$. We recall that $\alpha$ is a non-negative integer by assumption. Hence, $w(x,t)$ 
is continuous at $x = t$ and $w(t,t)= 0$, and so 
\begin{eqnarray*}
0 &=& \frac{\mathrm{d}}{\mathrm{d}t} \left(\int_0^t\mathrm{d}y - \int_{\infty}^t\mathrm{d}y\right) w(y) p_n(y)^2 \\
&=& \left( \lim_{x \to t^{-}} - \lim_{x \to t^{+}}\right) w(x)p_n^2(x) - \left(\int_0^t\mathrm{d}y - \int_{\infty}^t\mathrm{d}y\right) \frac{\partial}{\partial t} \left( p_n(y)^2  w(y) \right)\\ 
&=& -\alpha\int_0^{\infty} \frac{w(y)p_n^2(y)}{y-t} \mathrm{d}y + \int_0^{\infty} 2p_n \frac{\partial p_n}{\partial t} w(y)\mathrm{d}y
\end{eqnarray*}
where we have used
\[
\frac{\partial w(x,t)}{\partial t} = -\frac{\alpha w(x,t)}{x-t}.
\]
Using 
\[
\frac{\partial p_n}{\partial t} = \frac{\gamma_n'}{\gamma_n} \gamma_n x^n + \mathrm{lower \: order \: terms}
\]
gives us that 
\[
R_n = 2\frac{\gamma_n'}{\gamma_n}
\]
and hence 
\[
\frac{2}{a_n}\frac{\mathrm{d}a_n }{\mathrm{d}y} = 2
\left( \frac{\gamma_{n-1}}{\gamma_n} \right)^{-1}
\frac{\mathrm{d}}{\mathrm{d}t}
\left( \frac{\gamma_{n-1}}{\gamma_n} \right) = 2\frac{\gamma_{n-1}'}{\gamma_{n-1}} - 2\frac{\gamma_n'}{\gamma_n} = R_{n-1} - R_n.
\]

Similarly, differentiating \eqref{eq0:orthonormality} in the case of $i= j+1 = n$ with respect to $t$ shows
\begin{eqnarray*}
r_n &=& \int_0^{\infty} \left( \gamma_n' y^n + \gamma_{n,1}'y^{n-1} \right)p_{n-1}(y) \,\mathrm{d}y \\
&=& \frac{\gamma_n'}{\gamma_n} \int_0^{\infty} \gamma_n y^n p_{n-1} w(y)\mathrm{d}y +  \frac{\gamma_n'}{\gamma_n} \int_0^{\infty} \gamma_{n,1}'y^{n-1} p_{n-1} w(y)\mathrm{d}y.
\end{eqnarray*}
We use the fact that
\[
\gamma_n x^n = p_n(x) - \gamma_{n,1}x^{n-1} + \mathrm{lower \: order \: terms}
\]
in this expression to obtain 
\begin{eqnarray*}
r_n &=& \int_0^{\infty} \left( \gamma_{n,1}' - \frac{\gamma_{n,1}\gamma_n'}{\gamma_n}\right) y^{n-1} p_{n-1}w \, \mathrm{d}y\\
&=& \frac{\gamma_n}{\gamma_{n-1}} \int_0^{\infty} \left( \frac{\gamma_{n,1}'\gamma_n - \gamma_n' \gamma_{n,1}}{\gamma_n^2}\right) p_{n-1}^2 w \, \mathrm{d}y \\
&=& \frac{\mathrm{d}}{\mathrm{d}t} \left( \frac{\gamma_{n,1}}{\gamma_{n,1}}\right) 
\end{eqnarray*}
since
\[
b_n = \frac{\gamma_{n,1}}{\gamma_n} - \frac{\gamma_{n+1,1}}{\gamma_{n+1}} 
\]
the equation (\ref{eq3:Ldbn}) follows. 
\end{proof}

\begin{thm}\label{t3.8}
The coefficients, $R_n$ and $r_n$, satisfy the system of differential equations
\begin{subequations}
\begin{eqnarray}
\label{eq3:dRn}tR_n' &=& 2 r_n - \alpha + R_n( tR_n + 2n + \alpha + \mu - t)\\
\label{eq3:drn}t r_n' &=& \left( \frac{1- 2R_n}{R_n(1-R_n)} \right) r_n^2 - n(n+\alpha)\frac{R_n}{1-R_n}\\
&&  + (2n+\alpha+\mu)r_n + \frac{(2n+ \mu)r_n}{R_n-1} - \frac{\alpha r_n}{R_n}\nonumber 
\end{eqnarray}
\end{subequations}
\end{thm}

\begin{proof}
To obtain the first equation, note that
\[
\frac{\mathrm{d}b_n}{\mathrm{d}t} = R_n + t R_n' = r_n - r_{n+1}
\]
which gives \eqref{eq3:dRn} under the substitution of $r_{n+1}$ in accordance with \eqref{eq3:laddrec1}.

To obtain \eqref{eq3:drn} differentiate \eqref{eq3:anintermsofRnandrn}. Knowing $a_n'$ in terms of $a_n$, $r_n$ and $R_n$ from \eqref{eq3:Ldan} removes $a_n'$ while we may use \eqref{eq3:rnintermsofRnRnm1} to remove the remaining instances of $a_n^2$. Using \eqref{eq3:dRn} and \eqref{eq3:ladderc2} are used to eliminate $R_n'$ and $R_{n-1}$ to obtain an equivalent reformulation of \eqref{eq3:drn}.
\end{proof}

We observe that the differential equations
of corollary \ref{c3.4} and theorem \ref{t3.8} are identical upon setting
\begin{eqnarray*}
R_n &=& \frac{t+\theta_n}{t}\\
r_n &=& \frac{\kappa_n}{t} - \left(n+\frac{\mu}{2}\right).
\end{eqnarray*}
This demonstrates an equivalence between the characterization of the polynomial system
corresponding to (\ref{eq1:weight}) as implied by the method of isomonodromic deformation,
and the method of ladder operators. 

\subsection{Main results}

Recall that our remaining task is to relate the differential equations
of corollary to the Hamilton equations for a Painlev\'e V system.
The work in \cite{Forrester} and
\cite{BasorChen} both provide clues regarding relevant transformations. Explicitly, they suggest two M\"obius transforms
\begin{eqnarray*}
y &\sim& \frac{\theta_n}{t + \theta_n}\\
y &\sim& \frac{\theta_n + t}{\theta_n}
\end{eqnarray*}
both of which send $\infty$ to $1$, and send $0$ and $-t$ to $0$ and $\infty$ in different ways. 

\begin{proof}[Proof of proposition \ref{prop:maintheorem}]
We are required to show that $q$ satisfies \eqref{1.10a} using \eqref{eq2:forrestersthetarel}. We first find the derivatives of $q$ in terms of $\theta_n$ and $\kappa_n$,
\begin{eqnarray*}
q &=& \frac{\theta_n + t}{\theta_n}\\
q' &=& \frac{\theta_n - t\theta_n'}{\theta_n^2} \\
&=& - \frac{2\kappa_n + \theta_n(\theta_n + 2n + t + \alpha + \mu)}{\theta_n^2}\\
q'' &=& \frac{\theta_n'((2n+t+\alpha+\mu)\theta_n + 4\kappa_n)}{\theta_n^3} - \frac{1}{\theta_n} - 
\frac{2\kappa_n + \theta_n}{\theta_n^2}.
\end{eqnarray*}
However, using \eqref{eq2:forrestersthetarela} we have
\[
2\kappa_n = t\theta_n' - \theta_n(\theta_n+ 1 + 2n+ t + \alpha + \mu)
\]
giving
\begin{eqnarray*}
q'' &=& \frac{t(2t+2\theta_n)\theta_n'^2}{2\theta_n^3(t+\theta_n)} - \frac{(2t+\theta_n)\theta_n'}{\theta_n^2(t+\theta_n)} - \frac{1}{2\theta_n(t+\theta_n)} - \frac{\theta_n^2}{t(t+\theta_n)} \\
&+& \frac{\mu^2(t+\theta_n)^2 - \theta_n^2(t^2 + \alpha^2)}{2\theta_n^3(t+\theta_n)} + \frac{1}{2}\left( \frac{t}{t+\theta_n} - 2(1+2n+\alpha+\mu)\frac{t+\theta_n}{t\theta_n} - 5\right).
\end{eqnarray*}
Inverting the expression for $q$ in terms of $\theta_n$ gives us
\begin{eqnarray*}
\theta_n &=& \frac{t}{q-1}\\
\theta_n' &=& \frac{q - tq'-1}{(q-1)^2},
\end{eqnarray*}
and using these expressions show 
\begin{eqnarray*}
q'' &=& \left(\frac{1}{q-1} + \frac{1}{2q}\right)q'^2 - \frac{q'}{t} + \frac{(q-1)^2}{t^2}\left(\frac{\mu^2q}{2}-  \frac{\alpha^2}{2q} \right)\\
&&  - \frac{(1+2n+\alpha+\mu)q}{t} - \frac{q(1+q)}{2(q-1)},
\end{eqnarray*}
This is \eqref{1.10a} where
\begin{eqnarray*}
\alpha_1 = \frac{\mu^2}{2} \hspace{3cm} \alpha_2 = - \frac{\alpha^2}{2} \\
 \alpha_3 = - (2n +1 + \alpha + \mu) \hspace{2cm} \alpha_4 = -\frac{1}{2}.
\end{eqnarray*}
To obtain the corresponding $p$ variable, we remark that the equation for $q'$, as specified by the Hamiltonian in \eqref{eq1:HV}, is linear in $p$, hence determines $p$ uniquely in terms of $\theta_n$ and $\kappa_n$.
\end{proof}

For the differential equations to uniquely characterize $\theta_n, \kappa_n$, boundary values must be
specified. For this purpose, we note from (\ref{eq2:momentsf}) that the small $t$
leading order asymptotics for $\mu_k$ are
\begin{eqnarray*}
\mu_k &= & C_1 \left(1 - \frac{\mu}{k+\alpha+\mu}t +\frac{\mu(\mu-1)}{2(k+\alpha+\mu)(k+\alpha+\mu-1)}t^2 \right. \\
&& \left. -\frac{\mu(\mu-1)(\mu-2)}{6(k+\alpha+\mu)(k+\alpha+\mu-1)(k+\alpha+\mu-2)}t^3 + \ldots \right) \\
&& \hspace{.1cm} + C_2 t^{1+k+\alpha+\mu} \left(1- \frac{1+k+\alpha}{2+k+\alpha+\mu}t + \frac{(1+k+\alpha)(2+k+\alpha)}{2(2+k+\alpha+\mu)(3+k+\alpha+\mu)}t^2 \right. \\
&& \left. - \frac{(1+k+\alpha)(2+k+\alpha)(3+k+\alpha)}{6(2+k+\alpha+\mu)
(3+k+\alpha+\mu)(4+k+\alpha+\mu)}t^3 + \ldots \right).
\end{eqnarray*}
From this
the determinant (\ref{eq0:deltadetdef})
that defines $\Delta_n$ may be evaluated to leading orders by using the identity (see e.g.~\cite{No04}) 
\begin{equation}\label{eq3:detgamma}
\det(\Gamma(z_k+j))_{j,k=0,\ldots, n-1} = \prod_{k=0}^{n-1}\Gamma(z_k) \prod_{0 \leq j < k < n} (z_k-z_j).
\end{equation}
In particular, by letting $z_k = 1+\alpha+ \mu + k$ we have
\[
\Delta_n(0) = (1-\zeta)^n \prod_{k=1}^{n-1} k! \prod_{k=0}^{n-1} \Gamma(1+\alpha+\mu + k).
\]
Recalling (\ref{an}) then gives
\[
a_n^2(0) = n(n+\alpha+\mu),
\]
and knowing this \eqref{eq2:abdifferential} at $t=0$ implies
\[
b_n(0) = 2n+\alpha+1+\mu.
\]

To determine the rest of the expansion of \eqref{eq0:deltadetdef}, we first make use of
\eqref{eq3:detgamma} to compute the leading form of the analytic and non-analytic components as
\[
\Delta_n(t) = \Delta_n(0) \left( 1+ \frac{\mu}{\alpha+\mu}t +  O(t^2) + \chi_n t^{1+\alpha+\mu}( 1 + O(t) + O(t^{1+\alpha+\mu})) \right)
\]
where
\[
\chi_n = (1-\zeta)^{n-1}\left( (-1)^\mu + \frac{(\zeta-1)\sin(\pi \alpha)}{\sin(\pi(\alpha+\mu))}\right) \frac{ \Gamma(\mu+1) \Gamma(\alpha+\mu+n+1) }{(n-1)!\Gamma(\alpha+\mu + 1)\Gamma(\alpha+\mu+2)^2}.
\]
It follows from this and (\ref{an}) that the expansion of $a_n^2$ is
\begin{eqnarray*}
a_n^2 &=& n(n+\alpha+\mu) - \frac{n\mu^2(\alpha+\mu+n)}{(\mu+\alpha)^2}t^2 + O(t^3) \\
&& + t^{1+\alpha+\mu} \frac{\Gamma(\alpha+\mu+n+1)\Gamma(1+\mu)}{n!\Gamma(\alpha+\mu)\Gamma(\alpha+\mu+1)\Gamma(\alpha+\mu+2)} +  O(t^{1+\mu+\alpha}),
\end{eqnarray*}
and subsequently consistency with regards to \eqref{eq2:abdifferential} demands 
\begin{eqnarray*}
b_n &=& 2n + \alpha + \mu + 1  - \frac{\alpha}{\alpha+\mu}t + \frac{\alpha\mu(2n+\alpha+ \mu+1)}{(\alpha+\mu-1)(\alpha+\mu)^2(\alpha+1)}t^2 + O(t^3) \\
&&-t^{1+\alpha+\mu}\frac{\Gamma(\alpha+\mu+n+1)\Gamma(\mu+1)}{\Gamma(\alpha+\mu+1)^3 }+  O(t^{2+\mu+\alpha}).
\end{eqnarray*}
Substitution into (\ref{1.15}) then shows
\begin{eqnarray}
\theta_n(t) &=& - \frac{\mu}{\alpha+\mu}t + \frac{\alpha\mu(1+2n+\alpha+ \mu)}{(\alpha+\mu-1)(\alpha+\mu)^2(\alpha+1)}t^2 + O(t^3) \nonumber \\
&& \hspace{0.3cm} - t^{1+\alpha+\mu}\frac{\Gamma(\alpha+\mu+n+1)\Gamma(\mu+1)}{\Gamma(\alpha+\mu+1)^3 }+  O(t^{2+\mu+\alpha}) \nonumber \\
\kappa_n(t) &=& \frac{\mu(2n+\alpha+\mu)}{2+\alpha+\mu}t - \frac{2n\alpha\mu(n+\alpha+ \mu)}{((\alpha+\mu)^2-1)(\alpha+\mu)^2}t^2 + O(t^3) \nonumber \\
&& \hspace{.3cm}+ t^{1+\alpha+\mu} \frac{(\alpha+\mu+n-1)\Gamma(\mu+1)\Gamma(\alpha+\mu+n+1)}{\Gamma(\alpha+\mu+1)^3(\alpha + \mu)} \nonumber \\
&& +  O(t^{2+\mu+\alpha}) \label{3.19}
\end{eqnarray}
In particular, the small $t$ asymptotics for $q$ is therefore 
\begin{eqnarray*}
q = \frac{\theta_n+t}{\theta_n} &=&  -\frac{\alpha}{\mu} - \frac{\alpha(2n+\alpha+\mu+1)}{\mu(\alpha+\mu+1)(\alpha+\mu-1)}t + O(t^2)\\
&+& t^{\alpha+\mu} \frac{\Gamma(u)\Gamma(n+\alpha+\mu+1)}{\mu(\mu+\alpha)\Gamma(\alpha+\mu)^3}
+  O(t^{1+\mu+\alpha}).
\end{eqnarray*}

\end{document}